\documentclass[11pt,twoside, final]{amsart}
\usepackage{etoolbox,lastpage}
\usepackage{amsmath,amsthm,amscd,amsfonts,amssymb,enumerate}
\usepackage{graphicx}
\usepackage{color}
\usepackage[colorlinks]{hyperref}

\copyrightinfo{0}{Iranian Mathematical Society}
\pagespan{1}{\pageref*{LastPage}}
\commby{}
\newtheorem{theorem}{Theorem}[section]
\newtheorem{proposition}[theorem]{Proposition}
\newtheorem{lemma}[theorem]{Lemma}
\newtheorem{corollary}[theorem]{Corollary}
\theoremstyle{definition}
\newtheorem{definition}[theorem]{Definition}
\newtheorem{example}[theorem]{Example}
\theoremstyle{remark}

\numberwithin{equation}{section}

\input{tcilatex}

\begin{document}
\date{{\scriptsize Received: , Accepted: .}}
\title[On $(k,n)$-closed submodules]{On $(k,n)$-closed submodules}
\author[E. Yetkin Celikel]{Ece Yetkin Celikel}
\address[Ece Yetkin Celikel]{ Department of Mathematics, Faculty of Art and
Science, Gaziantep Unversity, Gaziantep, Turkey}
\email{yetkin@gantep.edu.tr}
\maketitle

\begin{abstract}
Let $R$ be a commutative ring with $1\neq 0$ and $M$ be an $R$-module. We
will call a proper submodule $N$ of $M$ as a semi $n$-absorbing submodule of 
$M$ if whenever $r\in R,$ $m\in M$ with $\ r^{n}m\in N$, then $r^{n}\in
(N:_{R}M)$ or $r^{n-1}m\in N$ . We will say $N$ to be a $(k,n)$-closed
submodule of $M$ if whenever $r\in R,$ $m\in M$ with $\ r^{k}m\in N$, then $%
r^{n}\in (N:_{R}M)$ or $r^{n-1}m\in N$. In this paper we introduce semi $n$%
-absorbing and $(k,n)$-closed submodules of modules over commutative rings,
and investigate their basic properties.

\end{abstract}

\section{\textbf{Introduction}}

Let $R$ be a commutative ring with $1\neq 0$ and $I$ be a proper ideal of $R$%
. As stated in \cite{andbad}, $I$ is called an $n$-absorbing (resp. strongly 
$n$-absorbing) ideal if whenever $x_{1}...x_{n+1}\in I$ for $%
x_{1},...,x_{n+1}\in R$ (resp. $I_{1}...I_{n+1}\subseteq I$ for ideals $%
I_{1},...,I_{n+1}$ of $R)$, then there are $n$ of the $x_{i}$'s (resp. $n$
of the $I_{i}$'s) whose product is in $I$. Recall thar a proper ideal $I$ of 
$R$ is said to be semi-prime ideal if whenever $r^{2}\in I$ for some $r\in R$%
, then $r\in I$. For generalizations of semi-prime ideals the reader may
consult \cite{Mo}. In \cite{andbad2}, D. F. Anderson and A. Badawi said $I$
to be a semi $n$-absorbing ideal if $x^{n+1}\in I$ for $x\in R$ implies $%
x^{n}\in I$. Also A. Badawi said that a proper ideal $I$ of $R$ is a $(m,n)$%
-closed ideal if $x^{m}\in I$ for $x\in R$ implies that $x^{n}\in I$ \cite%
{andbad2}. Let $M$ be an $R$-module. A proper submodule $N$ of $M$ is called 
$n$-absorbing (resp. strongly $n$-absorbing) submodule of $M$ if whenever $%
a_{1}...a_{n}m\in N$ for $a_{1},...,a_{n}\in R$ and $m\in M$ (resp. $%
I_{1}...I_{n}L\subseteq I$ for ideals $I_{1},...,I_{n}$ of $R$ and a
submodule $L$ of $M)$, then either $a_{1}...a_{n}\in (N:_{R}M)$ (resp. $%
I_{1}...I_{n}\subseteq (N:_{R}M))$ or there are $n-1$ of $a_{i}$'s ($I_{i}$%
's) whose product with $m$ (resp. $M)$ is in $N$ \cite{darso}. A proper
submodule $N$ of an $R$-module $M$ is called semi-prime if whenever $r\in R$
and $m\in M$ with $r^{2}m\in N$, then $rm\in N$. A proper submodule $N$ of $%
M $ is called a quasi-prime submodule of $M$ if whenever $a,b\in R$, $m\in M$
with $abm\in N$, then $am\in N$ or $bm\in N.$ More generally, we define $%
(k,n)$-closed submodules of an $R$-module $M$ as following: Let $R$ be a
commutative ring with identity and $k$, $n$ be positive integers. We call a
proper submodule $N$ of $M$ as a $(k,n)$-closed submodule of $M$ if whenever 
$r\in R,$ $m\in M$ with $\ r^{k}m\in N$, then $r^{n}\in (N:_{R}M)$ or $%
r^{n-1}m\in N$. In particular, we call $N$ as a semi $n$-absorbing submodule
of $M$ if whenever $r\in R,$ $m\in M$ with $\ r^{n}m\in N$, then $r^{n}\in
(N:_{R}M)$ or $r^{n-1}m\in N$. It is clear that a semi $n$-absorbing
submodule is $(n,n)$-closed.

Throughout we assume that all rings are commutative with $1\neq 0$, all
modules are considered to be unitary and $k,n$ are positive integers. The
radical of an ideal $I$ of $R$ is denoted by $\sqrt{I}.$ We denote the set
of invertible (unit) elements of $R$ by $U(R)$, i.e. $U(R)=\{u\in R:$ there
is a $v\in R$ such that $uv=vu=1_{R}\}.$ Let $N$ be a submodule of an $R$%
-module $M$. We will denote by $(N:_{R}M)$ \textit{the residual of $N$ by} $%
M $, that is, the set of all $r\in R$ such that $rM\subseteq N$. An $R$%
-module $M$ is called a \textit{multiplication module} if every submodule $N$
of $M$ has the form of $IM$ for some ideal $I$ of $R$. Note that, since $%
I\subseteq (N:_{R}M)$ then $N=IM\subseteq (N:_{R}M)M\subseteq N$. So that $%
N=(N:_{R}M)M$ \cite{ES}. For a submodule $N$ of $M$, if $N=IM$ for some
ideal $I$ of $R$, then we say that $I$ is a presentation ideal of $N$.
Clearly, every submodule of $M$ has a presentation ideal if and only if $M$
is a multiplication module. Let $N$ and $K$ be submodules of a
multiplication $R$-module $M$ with $N=I_{1}M$ and $K=I_{2}M$ for some ideals 
$I_{1}$ and $I_{2} $ of $R$. The product of $N$ and $K$ denoted by $NK$ is
defined by $NK=I_{1}I_{2}M$. Then by \cite[Theorem 3.4]{Am}, the product of $%
N$ and $K$ is independent of presentations of $N$ and $K$. Moreover, for $%
a,b\in M$, by $ab$, we mean the product of $Ra$ and $Rb$. Clearly, $NK$ is a
submodule of $M$ and $NK\subseteq N\cap K$ (see \cite{Am}). It is well-known
that if $R$ is a commutative ring and $M$ a non-zero multiplication $R$%
-module, then every proper submodule of $M$ is contained in a maximal
submodule of $M.$ \cite[Theorem 2.5]{ES}. As a generalization of Jacobson
radical of $R,$ the radical of the module $M$ is defined by the intersection
of all maximal submodules of $M$, that is $Rad(M)=\cap \{N:N$ is a maximal
submodule of $M\}.$ Let $N$ be a proper submodule of a non-zero $R$-module $M
$. Then the $M $-radical of $N$ denoted by $M$-$\mathrm{rad}(N)$ is defined
to be the intersection of all prime submodules of $M$ containing $N$. If $M$
has no prime submodule containing $N$, then we say $M$-$\mathrm{rad}(N)=M$.

In this study, we give many properties of $(k,n)$-closed submodules and also
obtain relationships among semi $n$-absorbing submodules, $(k,n)$-closed
submodules and the other concepts. For general background and terminology,
the reader may consult \cite{dd1} and \cite{S}.

\section{Properties of $(k,n)$-closed submodules}

In this section, we introduce and study basic properties of semi $n$%
-absorbing and $(k,n)$-closed submodules with many examples.

\begin{lemma}
\label{l1}Let $N$ be a proper submodule of an $R$-module $M.$ Then the
following statements are equivalent:
\end{lemma}

\begin{enumerate}
\item $N$ is a $(k,n)$-closed submodule of $M.$

\item If whenever $r\in R$ and $L$ is a submodule of $M$ with $%
r^{k}L\subseteq N$, then $r^{n-1}L\subseteq N$ or$\ r^{n}\in (N:_{R}M).$
\end{enumerate}

In particular, a proper submodule $N$ of $M$ is a semi $n$-absorbing
submodule of $M$ if and only if whenever $r\in R$, $L$ a submodule of $M$
with $r^{n}L\subseteq N$ implies either $r^{n}\in (N:_{R}M)$ or $%
r^{n-1}L\subseteq N$.

\begin{proof}
(1)$\implies $(2) Suppose that $N$ is a $(k,n)$-closed submodule of $M.$ Let 
$r\in R$ and $L$ be a submodule of $M$ with $r^{k}L\subseteq N.$ Assume that 
$r^{n-1}L\not\subseteq N.$ So $r^{n-1}m\notin N$ for some $m\in L.$ Since $%
r^{k}m\in N$ and $r^{n-1}m\notin N$, we conclude $r^{n}\in (N:_{R}M),$ as
needed.

(2)$\implies $(1) This part is clear.
\end{proof}

There are some relationships between $(k,n)$-closed submodules of $M$ and $%
(k,n)$-closed ideals of $R.$

\begin{theorem}
\label{t0}Let $M$ be an $R$-module, and $N$ be a proper submodule of $M.$ If 
$N$ is a $(k,n)$-closed submodule of $M,$ then $(N:_{R}M)$ is a $(k,n)$%
-closed submodule of $R.$ If $M$ is a multiplication $R$-module, then
presentation ideal of a $(k,n)$-closed submodule of $M$ is a $(k,n)$-closed
ideal of $R$.
\end{theorem}

\begin{proof}
Assume that $r\in R$ with $r^{k}\in (N:_{R}M)$ but $r^{n}\notin (N:_{R}M).$
Then there is an element $m\in M$ with $r^{n}m\notin N$ which means that $%
r^{n-1}m\notin N$. Since $r^{k}m\in N$, $r^{n-1}m\notin N$ and$\ r^{n}\notin
(N:_{R}M)$, this situation contradicts with our hypothesis. Thus $(N:_{R}M)$
is a $(k,n)$-closed ideal of $R.$
\end{proof}

However the converse of Theorem \ref{t0} is not true in general. For example
consider $N=6%
\mathbb{Z}
$ as a submodule of $%
\mathbb{Z}
$-module $%
\mathbb{Z}
.$ While $(N:_{%
\mathbb{Z}
}%
\mathbb{Z}
)=6%
\mathbb{Z}
$ is clearly a $(2,1)$-closed ideal of $%
\mathbb{Z}
,$ $N$ is not $(2,1)$-closed submodule of $%
\mathbb{Z}
.$ In fact $2^{2}.3^{2}\in N$ but $2^{1}\notin (N:_{%
\mathbb{Z}
}%
\mathbb{Z}
)=6%
\mathbb{Z}
$ and $2^{0}.3^{2}\notin N$.

\begin{theorem}
\label{tsm}Let $N$ be a proper submodule of $R$-module $M.$
\end{theorem}

\begin{enumerate}
\item If $N$ is a $(k,n)$-closed submodule of $M$, then $(N:_{R}m)$ is a $%
(k,n)$-closed ideal of $R$ for each $m\in M\backslash N.$

\item If $(N:_{R}m)$ is a $(k,n)$-closed ideal of $R$ for each $m\in
M\backslash N,$ then $N$ is a $(k,n+1)$-closed submodule of $M.$
\end{enumerate}

\begin{proof}
(1) Suppose that $r^{k}\in (N:_{R}m)$ and $r^{n}\notin (N:_{R}m)$ for some $%
m\in M\backslash N.$ Hence $r^{k}m\in N$ but $r^{n}m\notin N$ which means $%
r^{n-1}m\notin N.$ Since $N$ is a $(k,n)$-closed submodule of $M,$ we have $%
r^{n}\in (N:_{R}M)\subseteq (N:_{R}m),$ a contradiction. Thus $(N:_{R}m)$ is
a $(k,n)$-closed ideal of $R$ for each $m\in M\backslash N.$

(2) Let $r^{k}m\in N$ for $r\in R$ and $m\in M.$ Assume that $r^{n+1}\notin
(N:_{R}M)$. Since $r^{k}\in (N:_{R}m)$ and $(N:_{R}m)$ is a $(k,n)$-closed
ideal of $R$ for each $m\in M\backslash N$, we conclude that $r^{n}\in
(N:_{R}m).$ Therefore $r^{n}m\in N$. This means that $N$ is a $(k,n+1)$%
-closed submodule of $M.$
\end{proof}

\begin{lemma}
\label{lsm}Let $M$ be a finitely generated $R$-module such that $%
M=Rm_{1}+...+Rm_{t},$ $N$ be a proper submodule of $M$ and $k>n.$ Then
\end{lemma}

\begin{enumerate}
\item If $(N:_{R}m_{i})$ is a $(k,n)$-closed ideal of $R$ for all $%
i=1,...,t, $ then $(N:_{R}M)$ is a $(k,n)$-closed ideal of $R.$ In
particular, if $M=Rm$ be a cyclic $R$-module and $N$ is a proper submodule
of $M,$ then $(N:_{R}m)$ is a $(k,n)$-closed ideal of $R$ if and only if $%
(N:_{R}M)$ is a $(k,n)$-closed ideal of $R.$

\item Let $R$ be a division ring and $M=Rm$ be a cyclic $R$-module. Then $%
(N:_{R}m)$ is a $(k,n)$-closed ideal of $R$ if and only if $(N:_{R}m^{\prime
})$ is a $(k,n)$-closed ideal of $R$ for all elements $m^{\prime }\in M.$
\end{enumerate}

\begin{proof}
(1) Assume that $(N:_{R}m_{i})$ is a $(k,n)$-closed ideal of $R$ for all $%
i=1,...,t.$ Suppose that $r^{k}\in (N:_{R}M)$ and $r^{n}\notin (N:_{R}M)$
for some $r\in R$. Then $r^{n}\notin (N:_{R}m_{j})$ for some $j=1,...,t.$
Hence $r^{k}\notin (N:_{R}m_{j}),$ and so $r^{k}\notin (N:_{R}M)$ which
contradicts with our assumption. Thus $(N:_{R}M)$ is a $(k,n)$-closed ideal
of $R.$ The "in particular" part is clear.

(2) Suppose that $R$ is a division ring and $M=Rm$ is a cyclic $R$-module.
Then one can easily obtain that $(N:_{R}m)=(N:_{R}m^{\prime })$, so we are
done.
\end{proof}

\begin{theorem}
\label{c2}Let $R$ be a division ring and $N$ be a proper submodule of a
cyclic $R$-module $M=Rm.$
\end{theorem}

\begin{enumerate}
\item If $(N:_{R}m)$ is a $(k,n)$-closed ideal of $R,$ then $N$ is a $%
(k,n+1) $-closed submodule of $M.$

\item If $(N:_{R}m)$ is a semi $n$-absorbing ideal of $R,$ then $N$ is a
semi $(n+1)$-absorbing submodule of $M.$
\end{enumerate}

\begin{proof}
(1) From Theorem \ref{tsm} and Lemma \ref{lsm} (2), we are done.

(2) Since a semi $n$-absorbing ideal of $R$ is a $(n+1,n)$-closed ideal, $N$
is a $(n+1,n+1)$-closed submodule of $M$ by (1), so it is clear.
\end{proof}

In Theorem \ref{c2}, the condition "division ring" on $R$ is necessary.
Otherwise, if $(N:_{R}m)$ is a $(k,n)$-closed ideal of $R,$ then $N$ is not
need to be $(k,n+1)$-closed submodule of $M$ as in the following example.

\begin{example}
\label{e1}Consider $N=8%
\mathbb{Z}
$ as a submodule of $%
\mathbb{Z}
$-module $%
\mathbb{Z}
$. Then $(N:_{%
\mathbb{Z}
}1)=8%
\mathbb{Z}
$ is $(2,1)$-closed ideal but $N$ is not $(2,2)$-closed submodule of $M$. In
fact $2^{2}.2\in N$ but neither $2.2\in N$ nor $2^{2}\in (N:_{%
\mathbb{Z}
}%
\mathbb{Z}
).$
\end{example}

\begin{proposition}
Let $N$ be a proper submodule of an $R$-module $M$ and $k>t.$ Then the
following statements are equivalent:
\end{proposition}

\begin{enumerate}
\item $N$ is a $(k,n)$-closed submodule of $M$.

\item $(N:_{R}r^{k}m)=(N:_{R}r^{n-1}m)$ or $r^{n}\in (N:_{R}M)$ for $r\in R$
and $m\in M.$
\end{enumerate}

\begin{proof}
(1)$\Rightarrow $(2) Suppose that $N$ is a $(k,n)$-closed submodule of $M$
and $r^{n}\not\in (N:_{R}M).$ Let $s\in (N:_{R}r^{k}m).$ Hence $r^{k}(sm)\in
N$. Since $N$ is $(k,n)$-closed and $r^{n}\notin (N:_{R}M),$ we get $%
r^{n-1}sm\in N$. It follows $s\in (N:_{R}r^{n-1}m),$ that is $%
(N:_{R}r^{k}m)\subseteq (N:_{R}r^{n-1}m)$. Since the inverse inclusion is
always hold, this completes the proof.

(2)$\Rightarrow $(1) Suppose that $r\in R$, $m\in M$ with $\ r^{k}m\in N$.
If $r^{n}\in (N:_{R}M)$, then we are done. So assume that $\
(N:_{R}r^{k}m)=(N:_{R}r^{n-1}m).$ Thus $r^{n-1}m\in N$, as needed.
\end{proof}

The relations among the concepts of semi-prime, semi-$n$-absorbing,
quasi-prime, $n$ absorbing submodules and $(k,n)$-closed submodules are
given by the following theorem.

\begin{theorem}
\label{t1}Let $M$ be an $R$-module and $N$ be a proper submodule of $M.$
Then the following statements hold:
\end{theorem}

\begin{enumerate}
\item Let $N$ be a semi-prime submodule of $M$. Then $N$ is a $(k,n)$-closed
submodule of $M$ for all positive integers $k$ and $n.$ Moreover $N$ is a
semi $n$-absorbing submodule of $M$ for all positive integer $n.$

\item If $N$ is an $n$-absorbing submodule of $M$, then $N$ is a semi $n$%
-absorbing submodule of $M$ .

\item If $N$ is an $n$-absorbing submodule of $M$, then $N$ is a $(k,n)$%
-closed submodule of $M$ for every positive integer $k.$

\item If $N$ is a $(k,n)$-closed submodule of $M$, then $N$ is a $%
(k_{1},n_{1})$-closed submodule of $M$ for all $k_{1}\leq k$ and $n_{1}\geq
n.$

\item If $N$ is a semi $n$-absorbing submodule of $M,$ then $N$ is a semi $%
n_{1}$-absorbing submodule of $M$ for all $n_{1}\geq n.$

\item If $N$ is a quasi-prime submodule of $M$, then $N$ is a $(k,n)$-closed
submodule of $M$ for all positive integers $k\geq $ $n\geq 2.$
\end{enumerate}

\begin{proof}
(1), (2), (3) and (4) are clear from the definitions.

(5) Induction method on $n.$ For $n=1$, it is clear$.$ So suppose that $%
n\geq 2$ and $N$ is a semi $(n-1)$-absorbing submodule of $M.$ We show that $%
N$ is semi $n$-absorbing. Let $r\in R$ and $m\in M$ with $r^{n}m\in N$.
Assume that $r^{n}\in (N:_{R}M).$ Hence $r^{n-1}(rm)\in N$ which implies
that $r^{n-2}(rm)=r^{n-1}m\in N$ by introduction hypothesis. Thus $N$ is a
semi $n$-absorbing of $M$ for all $n_{1}\geq n.$

(6) We show that $N$ is a $(k,2)$-closed submodule of $M$ for all $k\geq 2$
by using mathematical induction on $k.$ Suppose that $N$ is a quasi-prime
submodule of $M$. Then $N$ is a $(k,2)$-closed submodule of $M$ for $k=2$
directly from their definitions. Now suppose that $N$ is a $(t,2)$-closed
submodule of $M$ for all $2\leq t<k$ and our aim is to show that $N$ is $%
(k,2)$-closed. Let $r^{k}m\in N$ for $r\in R$ and $m\in M.$ Assume that $%
r^{2}\notin (N:_{R}M).$ Since $r^{k-1}(rm)\in N,$ and $N$ is $(k-1,2)$%
-closed by induction hypothesis, we conclude that $r(rm)=r^{2}m\in N.$ Since 
$N$ is $(2,2)$-closed and $r^{2}\notin (N:_{R}M),$ we get $rm\in N.$ Thus $N$
is a $(k,2)$-closed submodule of $M$ for all $k\geq 2.$ Consequently, $N$ is
a $(k,n)$-closed submodule of $M$ for all positive integers $n$ with $k\geq
n\geq 2$ by (4).
\end{proof}

\begin{example}
\ The converses of (1)-(6) in Theorem \ref{t1} are not true in general as
these situations are shown in the following examples.

\begin{enumerate}
\item Let $N=30%
\mathbb{Z}
$ as a submodule of $%
\mathbb{Z}
$-module $%
\mathbb{Z}
.$ Since $N=2%
\mathbb{Z}
\cap 3%
\mathbb{Z}
\cap 5%
\mathbb{Z}
$ is intersection of semi-prime submodules of $%
\mathbb{Z}
,$ it is semi $2$-absorbing ($(2,2)$-closed) submodule of $%
\mathbb{Z}
$ from Theorem \ref{t2}. Also it is $(3,2)$-closed submodule of $%
\mathbb{Z}
$ from Theorem \ref{t1} (4). However $N$ is not $2$-absorbing submodule of $%
\mathbb{Z}
$. In fact $2.3.5\in N$ but $2.3\notin (N:_{%
\mathbb{Z}
}%
\mathbb{Z}
)$ and $2.5\notin N$ and $3.5\notin N$. So the converses of (2) and (3) are
not true.

\item Consider the submodule $N=\left( \overline{0}\right) $ of $%
\mathbb{Z}
$-module $%
\mathbb{Z}
_{p^{n}}$ where $p$ is a prime and $n$ is positive integer. Then $N$ is a $%
(n,n)$-closed submodule of $%
\mathbb{Z}
_{p^{n}}$, but $N$ is not $(n,n-1)$-closed as $p^{n}\overline{1}=\overline{0}%
\in N$ but neither $p^{n-2}\overline{1}\in N$ nor $p^{n-1}\in (N:_{%
\mathbb{Z}
}%
\mathbb{Z}
_{p^{n}})=(p^{n}).$ Note that $N$ in $%
\mathbb{Z}
$-module $%
\mathbb{Z}
_{p^{n}}$ is a semi $n$-absorbing submodule of $%
\mathbb{Z}
_{p^{n}}$, but it is not quasi-prime as $p^{n}\overline{1}\in N$ but $p%
\overline{1}\not\in N.$ Also it is not semi $(n-1)$-absorbing (it is also
not semi-prime clearly) submodule as $p^{n-1}p\in N$ but neither $p^{n-1}\in
(N:_{%
\mathbb{Z}
}%
\mathbb{Z}
)$ nor $p^{n-2}p=p^{n-1}\in N.$ Thus the coverses of (1), (4), (5) and (6)\
are not true.
\end{enumerate}
\end{example}

\begin{theorem}
\label{ti} Let $N$ be a proper submodule of $M.$ If $N$ is a semi $n$%
-absorbing submodule of $M$, then $N$ is a $(k,n)$-closed submodule of $M$
for all positive integer $k.$
\end{theorem}

\begin{proof}
If $k\leq n$, the the claim is clear. So suppose that $k>n$ and say $t:=k-n.$
Let $r^{k}m\in N$ for some $r\in R$ and $m\in M$. Assume that $r^{n}\notin
(N:_{R}M).$ Hence $r^{n}(r^{t}m)\in N$. Since $N$ is semi $n$-absorbing and $%
r^{n}\notin (N:_{R}M),$ we get $r^{n-1}(r^{t}m)=r^{n}(r^{t-1}m)\in N.$ This
follows $r^{n-1}(r^{t-1}m)=r^{n}(r^{t-2})\in N$ as again $N$ is a semi $n$%
-absorbing submodule of $M.$ It implies that $r^{n}(r^{t-3}m)\in N$. So we
continue with this argument and obtain that $r^{n}m\in N$ at the $t^{th}$
step. Finally we conclude $r^{n-1}m\in N$ which means that $N$ is a $(k,n)$%
-closed submodule of $M$.
\end{proof}

\begin{corollary}
\label{ciff}Let $N$ be a proper submodule of $M$ and $k>n$. Then $N$ is a $%
(k,n)$-closed submodule of $M$ if and only if $N$ is a semi $n$-absorbing
submodule of $M$.
\end{corollary}

\begin{proof}
Suppose that $N$ is $(k,n)$-closed and $r^{n}m\in N$ for $r\in R$ and $m\in
M.$ So $r^{k}m\in N$, and this implies that either $r^{n}\in (N:_{R}M)$ or $%
r^{n-1}m\in N.$ Thus $N$ is a semi $n$-absorbing submodule of $M.$ The
converse part follows from Theorem \ref{ti}.
\end{proof}

\begin{theorem}
\label{t2} Let $\{N_{\lambda }\}_{\lambda \in \Lambda }$ be a family of
semi-prime submodules of $M$. Then $\cap _{\lambda \in \Lambda }N_{\lambda }$
is a $(k,n)$-closed submodule of $M$ for all positive integers $k$ and $n$.
\end{theorem}

\begin{proof}
Suppose that $r^{k}m\in \cap _{\lambda \in \Lambda }N_{\lambda }$ for $r\in
R $ and $m\in M$. Then $r^{k}m\in N_{\lambda }$ for all $\lambda \in \Lambda
.$ Since each $N_{i}$ is semi-prime, we conclude that $rm\in N_{\lambda }$
for all $\lambda \in \Lambda .$ Thus $rm\in \cap _{\lambda \in \Lambda
}N_{\lambda }$ which means that $r^{n-1}m\in \cap _{\lambda \in \Lambda
}N_{\lambda }$ for all $n.$ From Theorem \ref{t1} (4), $\cap _{\lambda \in
\Lambda }N_{\lambda }$ is $(k,n)$-closed for all integers $k$ and $n.$
\end{proof}

\begin{corollary}
Let $N$ be a proper submodule of an $R$-module $M.$ Then $M-rad(N)$ and $%
Rad(M)$ are $(k,n)$-closed submodule of $M$ for all integers $k$ and $n.$
\end{corollary}

\begin{proof}
The result is clear from Theorem \ref{t2}.
\end{proof}

\begin{lemma}
\label{l3}\cite{Darat} Let $R$ be a commutative ring, $M$ a finitely
generated multiplication $R$-module and $N_{1},...,N_{t}$ are pairwise
comaximal $R$-submodules of $M$. Then the following statements hold:
\end{lemma}

\begin{enumerate}
\item $N_{1}N_{2}=N_{1}\cap N_{2}.$

\item $N_{1}\cap ...\cap N_{t-1}$ and $N_{t}$ are comaximal.

\item $N_{1}...N_{t}=N_{1}\cap ...\cap N_{t}.$
\end{enumerate}

\begin{theorem}
\label{tf2}Let $M$ is finitely generated multiplication $R$-module and $%
N_{1},...,N_{t}$ be semi-prime submodules of $M$. If $N_{1},...,N_{t}$ are
pairwise comaximal, then $N_{1}...N_{t}$ is a $(k,n)$-closed submodule of $M$
for all positive integers $k$ and $n$. In particular, if $N$ is semi-prime,
then $N^{n}$ is a $(k,n)$-closed submodule of $M.$
\end{theorem}

\begin{proof}
The claim is clear from Theorem \ref{t2} and Lemma \ref{l3}.
\end{proof}

D.F. Anderson and A. Badawi proved in Theorem 2.3 \cite{andbad2} that the
intersection of two semi $n$-absorbing ideals is also a semi $n$-absorbing
ideal of $R.$ However this situation is not true for submodules of any
module. The intersection of two semi $n$-absorbing submodule may not to be
semi $n$-absorbing as the following:

\begin{example}
Consider $%
\mathbb{Z}
$ as $%
\mathbb{Z}
$-module and two submodules $N=p^{n}%
\mathbb{Z}
$ and $K=q^{n}%
\mathbb{Z}
$ of $%
\mathbb{Z}
$ where $p$ and $q$ are prime integers. Clearly both of them are semi $n$%
-absorbing submodules of $%
\mathbb{Z}
.$ However $N\cap K=p^{n}q^{n}%
\mathbb{Z}
$ is not semi $n$-absorbing since $p^{n}(q^{n})\in N\cap K$ but $%
p^{n-1}(q^{n})\notin N\cap K$ and $p^{n}\notin (N\cap K:_{%
\mathbb{Z}
}%
\mathbb{Z}
).$
\end{example}

\begin{theorem}
Let $\{N_{\lambda }\}_{\lambda \in \Lambda }$ be a chain of $(k,n)$-closed
submodules of an $R$-module $M.$ Then $\tbigcap\limits_{\lambda \in \Lambda
}^{{}}N_{\lambda }$ is a $(k,n)$-closed submodule of $M.$
\end{theorem}

\begin{proof}
Let $r^{k}m\in N$ for $r\in R$ and $m\in M.$ If $r^{n}\in (N_{\lambda
}:_{R}M)$ for all $\lambda \in \Lambda $, then $r^{n}\in \cap (N_{\lambda
}:_{R}M)=(\cap N_{\lambda }:_{R}M)$, we are done. Supose that $r^{n}\notin
(N_{\lambda _{0}}:_{R}M)$ for some $\lambda _{0}\in \Lambda .$ Then $%
r^{n}\notin (N_{\lambda }:_{R}M)$ for all $N_{\lambda }\subseteq N_{\lambda
_{0}}.$ Hence $r^{n-1}m\in N_{\lambda }$ for all $N_{\lambda }\subseteq
N_{\lambda _{0}}$ as each $N_{\lambda }$ is $(k,n)$-closed. Therefore $%
r^{n-1}m\in \tbigcap\limits_{\lambda \in \Lambda }^{{}}N_{\lambda }$ which
means that $\tbigcap\limits_{\lambda \in \Lambda }^{{}}N_{\lambda }$ is a $%
(k,n)$-closed submodule of $M.$
\end{proof}

\begin{theorem}
Let $N_{1}$ and $N_{2}$ be proper submodules of an $R$-module $M.$
\end{theorem}

\begin{enumerate}
\item If $N_{1}$ is a semi $n_{1}$-absorbing and $N_{2}$ is a semi $n_{2}$%
-absorbing submodule of $M,$ then $N_{1}\cap N_{2}$ is semi $(n+1)$-closed
submodule of $M$ where $n=\max \{n_{1},n_{2}\}.$

\item If $N_{1},...,N_{t}$ be semi $n$-absorbing submodules of $M$. Then $%
N_{1}\cap ...\cap N_{t}$ is a semi $(n_{t}+t)$-absorbing submodule of $M$.

\item If $N_{1},...,N_{t}$ be semi $n_{t}$-absorbing submodules of $M$. Then 
$N_{1}\cap ...\cap N_{t}$ is a semi $(n_{t}+2)$-absorbing submodule of $M$
where $n=\max \{n_{1},...,n_{t}\}.$
\end{enumerate}

\begin{proof}
(1) Let $r\in R$ and $m\in M$ such that $r^{n+1}m\in N_{1}\cap N_{2}.$ First
observe from Corollary \ref{ciff} that $N_{1}$ and $N_{2}$ are $(n,n_{1})$%
-closed and $(n,n_{2})$-closed submodules of $M$, respectively. Hence we
have $r^{n_{1}}\in (N_{1}:_{R}M)$ or $r^{n_{1}-1}m\in N_{1}$ and $%
r^{n_{2}}\in (N_{2}:_{R}M)$ or $r^{n_{2}-1}m\in N_{2}.$ If $r^{n_{1}}\in
(N_{1}:_{R}M)$ and $r^{n_{2}}\in (N_{2}:_{R}M),$ then $r^{n}\in
(N_{1}:_{R}M)\cap (N_{2}:_{R}M)=(N_{1}\cap N_{2}:_{R}M)$. If $r^{n_{1}}\in
(N_{1}:_{R}M)$ and $r^{n_{2}-1}m\in N_{2},$ then $r^{n}m\in N_{1}\cap N_{2}.$
If symmetrically $r^{n_{1}-1}m\in N_{1},$ and $r^{n_{2}}\in (N_{2}:_{R}M)$,
then \ again we have $r^{n}m\in N_{1}\cap N_{2}.$ For the last, if $%
r^{n_{1}-1}m\in N_{1}$ and $r^{n_{2}-1}m\in N_{2},$ then $r^{n-1}m\in
N_{1}\cap N_{2}.$ Thus we conclude either $r^{n+1}\in (N_{1}\cap
N_{2}:_{R}M) $ or $r^{n}m\in (N_{1}\cap N_{2}:_{R}M),$ as needed.

(2) One can easily obtain the proof by using induction method on $t.$

(3) We use induction method on $t.$ If $t=3$, then the claim is clear from
(1) and (2). So assume that $t>3$ and the claim is satisfied for $t-1.$ Then 
$N_{1}\cap ...\cap N_{t-1}$ is semi $(n_{t-1}+2)$-absorbing. If $%
n_{t-1}+2<n_{t},$ then $N_{1}\cap ...\cap N_{t}$ is semi $(n_{t-1}+1)$%
-absorbing submodule of $M$ by part (1)$.$ Thus $N_{1}\cap ...\cap N_{t}$ is
semi $(n_{t-1}+2)$-absorbing submodule of $M$ by Theorem \ref{t1} (5). If $%
n_{t-1}+2=n_{t}$, then $N_{1}\cap ...\cap N_{t}$ is semi $(n_{t-1}+2)$%
-absorbing submodule of $M$ by part (2). If $n_{t-1}+2>n_{t},$ then $%
N_{1}\cap ...\cap N_{t}$ is $(n_{t-1}+3)$-absorbing by part (1). Here
observe that $n_{t-1}+3=n_{t}+2$ as $n_{t-1}+2>n_{t}$ and $n_{t-1}<n_{t}.$
Therefore $N_{1}\cap ...\cap N_{t}$ is semi $(n_{t}+2)$-absorbing submodule
of $M.$
\end{proof}

\begin{theorem}
Let $R$ be a division ring, $M$ be cyclic $R$-module, and $N_{1},...,N_{t}$
be $(k_{j},n_{j})$-closed submodules of $M$. Then $N_{1}\cap ...\cap N_{t}$
is a $(k,n+1)$-closed submodule of $M$ for all integers $k\leq \min
\{k_{1},...,k_{t}\}$ and $n\geq \min \{k,\max \{n_{1},...,n_{t}\}\}$.
\end{theorem}

\begin{proof}
Suppose that $N_{1},...,N_{t}$ are $(k_{j},n_{j})$-closed submodules of $M$.
Hence $(N_{1}:_{R}M),...,(N_{t}:_{R}M)$ are $(k_{j},n_{j})$-closed ideals of 
$R$ by Theorem \ref{t0}. Then $\tbigcap\limits_{j=1}^{t}(N_{j}:_{R}M)=(%
\tbigcap\limits_{j=1}^{t}N_{j}:_{R}M)$ is a $(k,n)$-closed ideal of $R$ for $%
k\leq \min \{k_{1},...,k_{t}\}$ and $n\geq \min \{k,\max
\{n_{1},...,n_{t}\}\}$ by Theorem 2.3 in \cite{andbad2}. Thus we conclude
that $\tbigcap\limits_{j=1}^{t}N_{j}$ is a $(k,n+1)$-closed submodule of $M$
by Theorem \ref{c2}.
\end{proof}

A non-zero submodule $N$ of an $R$-module $M$ is called a secondary
submodule of $M$ if for each $r\in R$ the homothety $N\overset{r}{%
\rightarrow }N$ is surjective or nilpotent (resp. surjective or zero). In
this case $P=\sqrt{(0:_{R}N)}$ is a prime ideal, and we call $N$ a $P$%
-secondary submodule of $M$. For more details concerning secondary submodule
of a module refer to \cite{mac}.

\begin{theorem}
\label{tsec}Let $N$ be a secondary submodule of an $R$-module $M.$ If $K$ is
a semi $n$-absorbing submodule of $M,$ then $N\cap K$ is a secondary
submodule of $M.$
\end{theorem}

\begin{proof}
Suppose that $N$ is a $P$-secondary submodule of $M$ and $r\in R.$ If $r\in
P=\sqrt{(0:_{R}N)}$, then clearly $r\in \sqrt{(0:_{R}N\cap K)}.$ So assume
that $r\notin P.$ Since $r^{n}\notin P,$ this implies that $r^{n}N=N.$ It is
needed to show that $r(N\cap K)=(N\cap K).$ Let $m\in N\cap K$. Since $%
N=r^{n}N,$ there is an element $m_{1}$ of $N$ such that $m=r^{n}m_{1}\in
N\cap K\subseteq K.$ Since $K$ is semi $n$-absorbing, we conclude either $%
r^{n}\in (K:_{R}M)$ or $r^{n-1}m\in K.$ If $r^{n}\in (K:_{R}M)$, then $%
N=r^{n}N\subseteq K$, and so $r(N\cap K)=rN=N\cap K.$ If $r^{n-1}m\in K$,
then $m=r^{n}m_{1}\in r(N\cap K),$ we are done.
\end{proof}

\begin{corollary}
\label{csec}Let $N$ and $K$ be proper submodules of an $R$-module $M$ with $%
K\subseteq N.$ If $N$ is a secondary semi $n$-absorbing submodule of $M,$
then $K$ is a semi $n$-absorbing submodule of $M.$
\end{corollary}

\begin{proof}
This is a direct consequence of Theorem \ref{tsec}.
\end{proof}

Let $N$ and $K$ be submodules of $M$ with $K\subseteq N.$ If $N$ is a semi $%
n $-absorbing submodule of $M$, then $K$ is not need to be a semi $n$%
-absorbing submodule of $M$ as the following example verifying this case. So
Example \ref{e} shows that the condition "secondary" in Corollary \ref{csec}
is necessary.

\begin{example}
\label{e}Consider a submodule $N=4%
\mathbb{Z}
$ of $%
\mathbb{Z}
$-module $%
\mathbb{Z}
$ and $K=12%
\mathbb{Z}
$. Then $K$ is clearly a semi $2$-absorbing submodule and $K\subseteq N,$
but $N$ is not semi $2$-absorbing submodule of $M$ as $2^{2}.3\in K$ but $%
2^{2}\notin (K:_{R}M)$ and $2.3\notin K.$
\end{example}

Let $R$ be an integral domain. Recall that if for every element $r$ of its
field of fractions $F$, at least one of $r$ or $r^{-1}$ belongs to $R$, then 
$R$ is called valuation domain.

\begin{proposition}
Let $R$ be a valuation domain with quotient field $K$. Let $M$ be an $R$%
-module and $N$ a proper submodule of $M$ Then $N$ is a semi $n$-absorbing
submodule of $M$ if and only if whenever $r\in K$, $m\in M$ with $r^{n+1}\in
N$ implies that $r^{n}m\in N.$
\end{proposition}

\begin{proof}
Suppose that $N$ is a semi $n$-absorbing submodule of $M$. Assume that $%
r^{n+1}m\in N$ but $r^{n+1}\notin (N:_{R}M)$ for some $r\in K$ and $m\in M.$
If $r\in R$, then we are done. So assume that $r\notin R.$ Since $R$ is a
valuation domain, $r^{-1}\in R.$ Hence we get $r^{-1}(r^{n+1}m)=r^{n}m\in N$%
. The converse part is clear.
\end{proof}

\begin{definition}
Let $N$ be a proper submodule of $M.$
\end{definition}

\begin{enumerate}
\item $N$ is said to be strongly semi $n$-absorbing submodule if whenever $I$
is an ideal and $L$ is a submodule of $M$ with $I^{n}L\subseteq N$ implies
that $I^{n}\subseteq (N:_{R}M)$ or $I^{n-1}L\subseteq N$.

\item $N$ is said to be strongly $(k,n)$-closed submodule if whenever $I$ is
an ideal and $L$ is a submodule of $M$ with $I^{k}L\subseteq N$ implies that 
$I^{n}\subseteq (N:_{R}M)$ or $I^{n-1}L\subseteq N.$
\end{enumerate}

Note that every strongly $(k,n)$-closed submodule is a $(k,n)$-closed
submodule of $M.$ Clearly a $(k,1)$-closed submodule is also a strongly $%
(k,1)$-closed submodule of $M.$ Also observe that a strongly semi $n$%
-absorbing submodule is a $n$-absorbing submodule of $M.$

\begin{lemma}
\label{l2}Let $N$ be a proper submodule of $M$. Then the following
statements are equivalent:
\end{lemma}

\begin{enumerate}
\item $N$ is a strongly $(k,n)$-closed submodule of $M.$

\item If $I$ is an ideal of $R$ and $m\in M$ with $I^{k}m\subseteq N$, then $%
I^{n}\subseteq (N:_{R}M)$ or $I^{n-1}m\subseteq N$.
\end{enumerate}

\begin{proof}
(1)$\implies $(2) It is obvious.

(2)$\implies $(1) Suppose that $I^{k}L\subseteq N$ for an ideal $I$ of $R$
and a submodule $L$ of $M$. Assume that $I^{n-1}L\nsubseteq N.$ Then there
is an element $m$ of $L$ such that $I^{n-1}m\nsubseteq N$ for some $m\in L.$
Since $I^{k}m\subseteq N$, we have $I^{n}\subseteq (N:_{R}M)$ by (2). Thus $%
N $ is a strongly $(k,n)$-closed submodule of $M.$
\end{proof}

\begin{theorem}
Let $R$ be a principal ideal domain and $N$ be a proper submodule of\ an $R$%
-module $M.$ Then the followings are equivalent:
\end{theorem}

\begin{enumerate}
\item $N$ is a $(k,n)$-closed submodule of $M.$

\item $N$ is a strongly $(k,n)$-closed submodule of $M.$
\end{enumerate}

\begin{proof}
(1)$\implies $(2) Since $I$ is principal, $I=(a)$ for some $a\in R.$\ So we
are done by Lemma \ref{l2}.

(2)$\implies $(1) It is clear.
\end{proof}

\begin{proposition}
Let $N$ be a proper submodule of an $R$-module $M.$ If $N$ is a $(k,n)$%
-closed submodule of $M$, then $(N:_{M}I)=\{m\in M:$ $Im\subseteq N\}$ is a $%
(k,n)$-closed submodule of $M$ for all ideal $I$ of $R.$ Moreover if $N$ is
a strongly $(k,n)$-closed submodule of $M$, then $%
(N:_{M}I^{k})=(N:_{M}I^{n-1}).$
\end{proposition}

\begin{proof}
Suppose that $r^{k}m\in (N:_{M}I)$ for $r\in R$ and $m\in M$. Hence $%
r^{k}Im\subseteq N$, which implies that either $r^{n}\in (N:_{R}M)$ or $%
r^{n-1}Im\subseteq N$ by Lemma \ref{l1}. This means $r^{n}\in
((N:_{R}M):_{R}I)=((N:_{M}I):_{R}M)$ or $r^{n-1}m\in (N:_{M}I).$ Thus $%
(N:_{M}I)$ is a $(k,n)$-closed submodule of $M$ for all ideal $I$ of $R.$
Now suppose that $N$ is a strongly $(k,n)$-closed submodule of $M.$ Since $%
(N:_{M}I^{n-1})\subseteq (N:_{M}I^{k})$ is always true, it is sufficient to
show the inverse inclusion. Let $m\in (N:_{M}I^{k})$. Then $I^{k}m\in N$,
and we have$\ I^{n}\subseteq (N:_{R}M)$ or $I^{n-1}m\in N$ from Lemma \ref%
{l2}$.$ If $I^{n-1}m\in N$, then $m\in (N:_{M}I^{n-1}),$ so we are done. So
suppose that $I^{n}\subseteq (N:_{R}M).$ Thus $I^{k}\subseteq (N:_{R}M),$ as
needed.
\end{proof}

\begin{theorem}
Let $N$ be a proper submodule of $M$. Then the following statements are
equivalent:
\end{theorem}

\begin{enumerate}
\item $N$ is a strongly $(k,n)$-closed submodule of $M.$

\item For any ideal $I$ of $R$ and $N\subseteq L$ a submodule of $M$ with $%
I^{k}L\subseteq N$ implies that $I^{n}\subseteq (N:_{R}M)$ or $%
I^{n-1}L\subseteq N.$
\end{enumerate}

\begin{proof}
(1)$\implies $(2) It is clear.

(2)$\implies $(1) Let $K$ be a submodule of $M$ and $I$ an ideal of $R$ such
that $I^{k}K\subseteq N.$ Hence $I^{k}(K+N)=I^{k}K+I^{k}N\subseteq N.$ Put $%
L=K+N.$ Since $N$ is strongly $(k,n)$-closed, we conclude that either $%
I^{n}\subseteq (N:_{R}M)$ or $I^{n-1}L\subseteq N$ by hypothesis (2). Thus $%
I^{n}\subseteq (N:_{R}M)$ or $I^{n-1}K\subseteq N.$
\end{proof}

\begin{theorem}
\label{st}Let $N$ be a $(k,2)$-closed submodule of $M,$ and $L$ a submodule
of $M.$ Then the following statements are satisfied:
\end{theorem}

\begin{enumerate}
\item If $L^{k}M\subseteq N$, then $2L^{2}\subseteq (N:_{R}M).$

\item If $2\in U(R),$ then $N$ is a strongly $(k,2)$-closed submodule of $M.$
\end{enumerate}

\begin{proof}
(1) Suppose that $L^{k}M\subseteq N$. Then $l_{1}^{k}m,$ $l_{2}^{k}m,$ $%
(l_{1}+l_{2})^{k}m\in N$ for all $m\in M$, for all $l_{1},l_{2}.$ Since $N$
is $(k,2)$-closed, we conclude that (either $l_{1}^{2}\in (N:_{R}M)$ or $%
l_{1}m\in N)$ and $($either $l_{2}^{2}\in (N:_{R}M)$ or $l_{2}m\in N)$ and
(either $(l_{1}+l_{2})^{2}\in (N:_{R}M)$ or $(l_{1}+l_{2})m\in N)$ which
means $l_{1}^{2}m,$ $l_{2}^{2}m,$ $(l_{1}+l_{2})^{2}\in N$. Then $%
2l_{1}l_{2}m=((l_{1}+l_{2})^{2}-l_{1}^{2}-l_{2}^{2})m\in N$. Thus $%
2L^{2}M\subseteq N,$ and so $2L^{2}\subseteq (N:_{R}M).$

(2) Let $2\in U(R).$ Since $2L^{2}M\subseteq N$ from (1), we conclude that $%
L^{2}\subseteq (N:_{R}M).$
\end{proof}

Now we extend well-known results about prime submodules, $n$-absorbing
submodules and $(m,n)$-closed ideals to $(k,n)$-closed submodules.

\begin{theorem}
\label{s-1}Let $N$ be a proper submodule of $M,$ and $S$ be a
multiplicatively closed subset of $R$ such that $(N:_{R}M)\cap S=\emptyset .$
If $N$ is a $(k,n)$-closed submodule of $M$, then $S^{-1}N$ is a $(k,n)$%
-closed submodule of $S^{-1}M.$ In particular, if $N$ is a semi $n$%
-absorbing submodule of $M$, then $S^{-1}N$ is a semi $n$-absorbing
submodule of $S^{-1}M.$
\end{theorem}

\begin{proof}
Let $\left( \frac{r}{s_{1}}\right) ^{k}\left( \frac{m}{s_{2}}\right) \in
S^{-1}N$. Hence $ur^{k}m\in N$ for some $u\in S$. Hence $(ur)^{k}m\in N$.
Since $N$ is $(k,n)$-closed, $(ur)^{n-1}m\in N$ or $(ur)^{n}\in (N_{:R}M)$
which follows either $\left( \frac{r}{s_{1}}\right) ^{n-1}\left( \frac{m}{%
s_{2}}\right) =\frac{u^{n-1}r^{n-1}}{u^{n-1}s_{1}^{n}}\frac{m}{s_{2}}\in
S^{-1}N$ or $\left( \frac{r}{s_{1}}\right) ^{n}=\frac{u^{n}r^{n}}{%
u^{n}s_{1}^{n}}\in S^{-1}(N:_{R}M)\subseteq (S^{-1}N:_{S^{-1}R}S^{-1}M)$.
"In particular" part is clear as a semi $n$-absorbing submodule is a $(n,n)$%
-closed submodule of $M$.
\end{proof}

\begin{corollary}
Let $S$ be a multiplicatively closed subset of$\ $such that $S\cap
(N:_{R}M)=\emptyset $ with $2\in S$. If $N$ is a strongly $(k,2)$-closed
submodule of $M$, then $S^{-1}N$ is a strongly $(k,2)$-closed submodule of $%
S^{-1}M$.
\end{corollary}

\begin{proof}
Let $S^{-1}K$ be a submodule of $S^{-1}M$ such that $(S^{-1}K)^{k}(S^{-1}M)%
\subseteq S^{-1}N.$ Since $2\in S$, $2\notin U(S^{-1}R),$ we are done by
Theorem \ref{st} (2).
\end{proof}

\begin{corollary}
Let $N$ be a proper submodule of $M$, and $P$ a prime submodule of $M$
containing $N$. Then $N$ is a $(k,n)$-closed submoule of $M$ if and only if $%
N_{P}$ is a $(k,n)$-closed submodule of $M_{P}.$
\end{corollary}

\begin{proof}
If $N$ is a $(k,n)$-closed submodule of $M$, then $N_{P}$ is a $(k,n)$%
-closed submodule of $M_{P}$ by Theorem \ref{s-1}. Conversely suppose that $%
r\in R$, $m\in M$ with $r^{k}m\in N$. Let $\Omega =\{u\in R:ur^{n}m\in N\}.$
Then $\left( \frac{r}{1}\right) ^{k}\frac{m}{1}\in N_{P}$ implies that $%
\left( \frac{r}{1}\right) ^{n-1}\frac{m}{1}\in N_{P}$ or $\left( \frac{r}{1}%
\right) ^{n}\in (N_{P}:_{R_{P}}M_{P})$ as $N_{P}$ is $(k,n)$-closed.
Therefore $ur^{n}m\in N_{P}$ for some $u\in R\backslash P.$ Hence $\Omega
\nsubseteq P.$ Also $\Omega \nsubseteq P^{\prime }$ where $P^{\prime }$ is
any prime submodule of $M$ with $I\nsubseteq P^{\prime }.$ Therefore $\Omega
=R$, which means that $r^{n}m\in N$. Thus $N$ is a $(k,n)$-closed submodule
of $M.$
\end{proof}

\begin{theorem}
Let $M$, $M^{\prime }$ be $R$-modules with unity, and $f:M\rightarrow
M^{\prime }$ an $R$-module homomorphism.
\end{theorem}

\begin{enumerate}
\item If $N^{\prime }$ is a $(k,n)$-closed (resp. semi $n$-absorbing)
submodule of $M^{\prime },$ then $f^{-1}(N^{\prime })$ is a $(k,n)$-closed
(resp. semi $n$-absorbing) submodule of $M$.

\item If $f$ \ is onto and $N$ is a $(k,n)$-closed (resp. semi $n$%
-absorbing) submodule of $M$ containing $Ker$ $f,$ then $f(N)$ is a $(k,n)$%
-closed (resp. semi $n$-absorbing) submodule of $M^{\prime }$
\end{enumerate}

\begin{proof}
The reader can easily obtain the proof, so\ it is omitted.
\end{proof}

\begin{corollary}
Let $M$, $M^{\prime }$ be $R$-modules and $N$, $K$ be proper submodules of $%
M.$ Then the following statements hold:
\end{corollary}

\begin{enumerate}
\item If $M\subseteq M^{\prime }$ and $N$ is a $(k,n)$-closed (resp. semi $n$%
-absorbing) submodule of $M^{\prime },$ then $N\cap M$ is a $(k,n)$-closed
(resp. semi $n$-absorbing) submodule of $M.$

\item If $K\subseteq N$, then $N/K$ is a $(k,n)$-closed (resp. semi $n$%
-absorbing) submodule of $M/K$ if and only if $K$ is a $(k,n)$-closed (resp.
semi $n$-absorbing) submodule of $M.$
\end{enumerate}

\begin{theorem}
Let $M_{1}$, $M_{2}$ be $R$-modules with $M=M_{1}\oplus M_{2},$ and let $%
N_{1}$, $N_{2}$ be proper submodules of $M_{1}$, $M_{2},$ respectively.
\end{theorem}

\begin{enumerate}
\item $N_{1}$ is a $(k_{1},n_{1})$-closed submodule of $M_{1}$ if and only
if $N_{1}\oplus M_{2}$ is a $(k,n)$-closed submodule of $M_{1}\oplus M_{2}$
for all positive integers $k_{1}\leq k$ and $n\geq n_{1}.$

\item $N_{2}$ is a $(k_{2},n_{2})$-closed submodule of $M_{2}$ if and only
if $M_{1}\oplus N_{2}$ is a $(k,n)$-closed submodule of $M_{1}\oplus M_{2}$
for all positive integers $k_{2}\leq k$ and $n\geq n_{2}.$
\end{enumerate}

\begin{proof}
(1) Suppose that $N_{1}$ is a $(k_{1},n_{1})$-closed submodule of $M_{1}$.
Assume that $r^{k_{1}}(m_{1},m_{2})\in N_{1}\oplus M_{2}$ but $%
r^{n_{1}-1}(m_{1},m_{2})\notin N_{1}\oplus M_{2}$. Then $r^{n_{1}-1}m_{1}%
\notin N_{1},$ which implies $r^{n_{1}}\in (N_{1}:_{R}M_{1}).$ Thus $%
r^{n_{1}}\in (N_{1}\oplus M_{2}:_{R}M).$ Consequently, $N_{1}\oplus M_{2}$
is a $(k,n)$-closed submodule of $M_{1}\oplus M_{2}$ for all positive
integers $k_{1}\leq k$ and $n\geq n_{1}$ by Theorem \ref{t1} (4). The
converse part can be obtained easily by using the similar argument.

(2) It can be easily verified similar to (1).
\end{proof}

\begin{theorem}
Let $M_{1}$, $M_{2}$ be $R$-modules, $N_{1}$ be a $(k_{1},n_{1})$-closed
submodule of $M_{1}$, and $N_{2}$ be a $(k_{2},n_{2})$-closed submodule of $%
M_{2}.$ Then $N_{1}\oplus N_{2}$ is a $(k,n)$-closed submodule of $%
M_{1}\oplus M_{2}$ for all positive integers $k\leq \min \{k_{1},k_{2}\}$
and $n\geq \max \{n_{1},n_{2}\}+1.$
\end{theorem}

\begin{proof}
Suppose that $r\in R$ and $(m_{1},m_{2})\in M$ such that $%
r^{k}(m_{1},m_{2})\in N_{1}\oplus N_{2}$. Hence $r^{k}m_{1}\in N_{1}.$ Since 
$r^{k_{1}}m_{1}\in N_{1}$ and $N_{1}$ is a $(k_{1},n_{1})$-closed submodule
of $M_{1},$ we have $r^{n_{1}}m_{1}\in N_{1}.$ Similarly, since $%
r^{k_{2}}m_{2}\in N_{2}$ and $N_{2}$ is a $(k_{2},n_{2})$-closed submodule
of $M_{2},$ we get $r^{n_{2}}m_{2}\in N_{2}.$ Thus $r^{n}m_{1}\in N_{1}$ and 
$r^{n}m_{2}\in N_{2}$ for all $n\geq \max \{n_{1},n_{2}\}.$ Therefore $%
r^{n}(m_{1},m_{2})\in N_{1}\oplus N_{2}$, as needed.
\end{proof}

D.F Anderson and A. Badawi determined in \cite{andbad2}.when the powers of
principal prime ideal or the ideals of the form $%
p_{1}^{t_{1}}...p_{i}^{t_{i}}$ where $p_{1},...,p_{t}$ are non associate
prime elements of $R$ and $t_{1},...,t_{n}$ are positive integers are $(m,n)$%
-closed ideal of an integral domain $R.$ Analogous to them, we may conclude
many results for submodules of multiplication modules over commutative
rings. Some of them are presented as the following.

\begin{theorem}
\label{tkn}Let $R$ be an integral domain and $M$ be a multiplication $R$%
-module. Let $(N:_{R}M)=p^{t}R$ where $p$ is prime element of $R$ and $k>0.$
If $N$ is a $(k,n)$-closed submodule of $M$, then the following statements
are hold:
\end{theorem}

\begin{enumerate}
\item $t=ka+r,$ where $a$ and $r$ are integers such that $a\geq 0,$ $1\leq
r\leq n$, $a(k$ $\func{mod}n)+r\leq n$, and if $a\neq 0$, then $k=n+c$ for
an integer $c$ with $1\leq c\leq n-1.$

\item If $k=bn+c$ for integers $b$ and $c$ with $b\geq 2$ and $0\leq c\leq
n-1$, then $t\in \{1,...,n\}.$ If $k=n+c$ for an integer $c$ with $0\leq
c\leq n-1$, then $t\in \cup _{h=1}^{n}\{ki+h:i\in 
\mathbb{Z}
$ and $0\leq ic\leq n-h\}.$
\end{enumerate}

\begin{proof}
Suppose that $N$ is a $(k,n)$-closed submodule of $M.$ Then $(N:_{R}M)$ is a 
$(k,n)$-closed ideal of $M$ by Theorem \ref{t0}. So we are done from Theorem
3.1 in \cite{andbad2}.
\end{proof}

\begin{corollary}
Let $M$ be a multiplication $R$-module where $R$ is an integral domain, and $%
(N:_{R}M)=p^{t}R$ where $p$ is prime element of $R$, $t>0.$ If $N$ is a semi 
$n$-absorbing submodule of $M$, then $t=na+r,$ where $a$ and $r$ are
integers such that $a\geq 0,$ $1\leq r<n$, that is $t\in \cup
_{h=1}^{n}\{ni+h:i\in 
\mathbb{Z}
$ and $0\leq i<n-h\}.$
\end{corollary}

\begin{proof}
Since a semi $n$-absorbing submodule is a $(n,n)$-closed submodule of $M,$
the resut is clear by Theorem \ref{tkn}.
\end{proof}

\begin{corollary}
Let $R$ be an integral domain and $(N:_{R}M)=p^{t}R$ where $p$ is a prime
element of $R$ and $t$ is a positive integer. Then $N$ is a semi 2-absorbing
submodule of $M$, then $t\in \{1,2\}.$
\end{corollary}

Consider a $%
\mathbb{Z}
$-module $M=%
\mathbb{Z}
$ and a submodule $N=2^{3}%
\mathbb{Z}
$ of $M$. It is shown in Example \ref{e1} that $N$ is not a semi 2-absorbing
submodule of $M$ for $t=3.$

\begin{theorem}
Let $R$ be a principal ideal domain, $N$ a proper submodule of a
multiplication $R$-module $M$ and $k$, $n$ be integers with $1\leq n\leq k.$
If $N$ is a (strongly) $(k,n)$-closed submodule of $M$, then $%
N=P_{1}^{t_{1}}...P_{i}^{t_{i}}$ where $P_{1},...,P_{i}$ are nonassociate
prime submodules of $M$, $t_{1},...,t_{i}$ are positive integers, and one of
the following two conditions holds:
\end{theorem}

\begin{enumerate}
\item If $k=bn+c$ for integers $b$ and $c$ with $b\geq 2$ and $0\leq c\leq
n-1$, then $t_{j}\in \{1,...,n\}$ for every $1\leq j\leq i.$

\item If $k=n+c$ for an integer $c$ with $0\leq c\leq n-1$, then $t_{j}\in
\cup _{h=1}^{n}\{kv+h:v\in 
\mathbb{Z}
$ and $0\leq vc\leq n-h\}$ for every $1\leq j\leq i.$
\end{enumerate}

\begin{proof}
Suppose that $N$ is $(k,n)$-closed submodule of $M.$ Then $(N:_{R}M)$ is a $%
(k,n)$-closed ideal of $R$ by Theorem \ref{t0}. Hence $%
(N:_{R}M)=p_{1}^{t_{1}}...p_{i}^{t_{i}}R$ for some nonassociate prime
elements of $R$, $t_{1},...,t_{i}$ are positive integers, and the conditions
(1) or (2) is satisfied for $k$ and $n$ by \cite{andbad2}. Thus $%
N=p_{1}^{t_{1}}...p_{i}^{t_{i}}M$. Put $P_{i}^{t_{i}}=p_{i}^{t_{i}}M$ for
all $i=1,...,t$, \ so $N=P_{1}^{t_{1}}...P_{i}^{t_{i}}$, we are done.
\end{proof}

\begin{theorem}
Let $N$ be a proper submodule of a multiplication $R$-module $M$ where $R$
is an integral domain and $k$, $n$ be integers with $1\leq n\leq k$. Suppose
that $N=P^{t}$, where $P$ is a prime submodule of $M$ and $t$ is a positive
integer. If $N$ is a $(k,n)$-closed submodule of $M,$ then one of the
following statements holds:
\end{theorem}

\begin{enumerate}
\item $1\leq t\leq n.$

\item There is a positive integer $a$ such that $t=ka+r=na+d$ for an integer 
$r$ and $d$ with $1\leq r$, $d\leq n-1.$

\item There is a positive integer $a$ such that $t=ka+r=n(a+1)$ for an
integer $r$ with $1\leq r\leq n-1.$
\end{enumerate}

\begin{proof}
From Theorem \ref{t0} and Theorem 3.8 in \cite{andbad2}, the result is clear.
\end{proof}

\end{document}